\documentclass[12pt,a4paper,draft]{article}
\usepackage[english]{babel}
\usepackage{cite}
\usepackage{latexsym,amsfonts,amssymb,amsmath,longtable,stmaryrd,mathrsfs,amsthm}
\usepackage[mathlines,pagewise]{lineno}
\renewcommand{\leq}{\leqslant}
\renewcommand{\geq}{\geqslant}

\renewcommand{\mathcal}{\mathscr}
\theoremstyle{plain}
\newtheorem{Lemma}{{\bfseries Lemma}}

\newtheorem{Cor}{{\bfseries Corollary}}
\newtheorem{Theo}{{\bfseries Theorem}}

 \DeclareMathOperator{\Hall}{Hall}

\newcommand{\F}{\mathbb F}

\newcommand{\E}{\mathcal E}

\newcommand{\D}{\mathcal D}

\DeclareMathOperator{\Alt}{Alt}
\DeclareMathOperator{\Sym}{Sym}
\newcommand{\arbitraryext}{\,\ldotp}
\newcommand{\nonsplitext}{\,{}^{\text{\normalsize{\textperiodcentered}}}}
\newcommand{\splitext}{\,\colon\!}

\title{\vspace{-1cm} 
Confirmation for Wielandt's conjecture\thanks{Research  is supported by
a NNSF grant of China (Grant \# 11371225)  and Wu Wen-Tsun Key Laboratory of Mathematics, USTC, Chinese Academy of Sciences.}}
\begin{document}

\author{Wenbin Guo\\
{\small Department of Mathematics, University of Science and
Technology of China,}\\ {\small Hefei 230026, P. R. China}\\
{\small E-mail:
wbguo@ustc.edu.cn}\\ \\
D. O. Revin, E. P. Vdovin\\
{\small Sobolev Institute of Mathematics and Novosibirsk State University,}\\
{\small Novosibirsk 630090, Russia}\\
{\small E-mail: revin@math.nsc.ru, vdovin@math.nsc.ru}}

 \date{}
\maketitle

\pagenumbering{arabic}
\begin{abstract}

Let $\pi$ be a set of primes. By H.Wielandt definition,  {\it Sylow  $\pi$-theorem} holds for a finite group  $G$ if  all maximal $\pi$-subgroups of  $G$  are conjugate. In the paper, the following statement is proven. Assume that $\pi$ is a union of disjoint subsets  $\sigma$ and $\tau$ and a finite group $G$ possesses a  $\pi$-Hall subgroup which is a direct product of a $\sigma$-subgroup and a $\tau$-subgroup. Furthermore, assume that both the Sylow $\sigma$-theorem and $\tau$-theorem hold for $G$. Then the  Sylow $\pi$-theorem holds for $G$. This result confirms a conjecture  posed by H.\,Wielandt in~1959.


\end{abstract}

\section*{Introduction}

In the paper, the  term  `group' always means  `finite group', $p$ is always supposed to be a prime, $\pi$ is a subset of the set of all primes, and $\pi'$ is its complement in the set of all primes.

Lagrange's theorem
states that the order $|G|$ of a group $G$ is divisible by the
order of every subgroup of $G$. This simple statement has extraordinary significance
and largely determines the problems of finite group theory. Lagrange's theorem
shows the extent to which the order of a group determines its subgroup structure.
For example, it turns out that every group of prime order is cyclic and
contains no non-trivial proper subgroups.

It is well-known that the converse of Lagrange theorem is
false. However,
Sylow's theorem states that if $p$ is a prime then
\begin{itemize}
\item every  group $G$ contains a so-called Sylow $p$-subgroup, i.\,e.  a subgroup  $H$ such that $|H|$ is a power of $p$ and the index $|G:H|$ is not divisible by $p$, and
\item every $p$-subgroup of $G$ (i.\,e. a subgroup whose order is a power of $p$) is conjugate with a subgroup of $H$.
\end{itemize}
Thus, the converse of Lagrange's theorem holds for certain
divisors of the group order. Moreover, it turns out that, in a finite group, the structure and properties
of every $p$-sub\-gro\-up are determined in many respects by the structure and
properties of any Sylow $p$-subgroup. In fact, Sylow's theorem is  considered by specialists as a cornerstone of finite group theory.

A natural generalization of the concept of a Sylow $p$-subgroup
is that of a $\pi$-Hall subgroup.

Recall that a subgroup $H$ of $G$ is called a {{\em $\pi$-Hall subgroup}} if
\begin{itemize}
\item all prime divisors of $|H|$ lie in $\pi$ (i.e. $H$ is a {\em $\pi$-subgroup}) and
\item all prime divisors of the index  $|G:H|$ of $H$ lie in $\pi'$.
\end{itemize}

According to P.Hall \cite{Hall}, we  say that a finite group $G$ {\em satisfies} $\D_\pi$ (or is a {\em $\D_\pi$-group}, or briefly $G\in\D_\pi$)
if  $G$ contains a  $\pi$-Hall subgroup $H$ and every $\pi$-subgroup of  $G$ is conjugate in $G$ with some subgroup of $H$. Thus, for a group $G$ the condition $\D_\pi$ means that the complete analog of Sylow's theorem holds for $\pi$-subgroups of $G$. That is why, together with Hall's notation, the terminology introduced by H.\,Wielandt in \cite{Wie3,Wie2} is used. According to Wielandt,  {\it the Sylow $\pi$-theorem} holds (der $\pi$-Sylow-Satz gilt) for a group $G$  if  $G$ satisfies $\D_\pi$. Sylow's theorem implies that the Sylow $\pi$-theorem for $G$ is equivalent to the conjugacy of all maximal $\pi$-subgroups of $G$.

Recall that a  group is nilpotent if and only if it is a direct product of its Sylow subgroups.
In \cite{Wie}, H.Wielandt proved the following theorem.

\begin{Theo}\label{W1954} {\rm \cite[Satz]{Wie}}
Assume that a  group $G$ possesses a nilpotent $\pi$-Hall subgroup for a set $\pi$ of primes. Then $G$ satisfies $\D_\pi$.
\end{Theo}

This result is regarded to be classical. It can be found in well-known textbooks \cite{Suz,ShemBook,DH,Isaacs,Robinson,GuoBook}. Wielandt mentioned the result
is his talk  at the XIII-th International Mathematical Congress in Edinburgh~\cite{Wie3}.

There are a lot of generalizations and analogs of Wielandt's theorem which was proved by many specialists
(see, for example, \cite{Baer,Hall,Hart,Zappa,Tibiletti,Shem3,Shem4,Rusakov,Rusakov1,Rusakov2,Gol1,Gol2,Tyshkevich,Wehrfritz}).
Wielandt's theorem plays an important role in the study of $\D_\pi$-groups (cf. \cite{MR,R,R2,R3,R4,R5,R6,VR1,VR2,VR3,VR4,Moreto,GuoSkiba,LiGuoKhan,Gross1,Gross2,Gross3,Gross4}).

One of the earliest generalizations of Wielandt's theorem was obtained by Wielandt himself in \cite{Wie2}:

\begin{Theo}\label{W1959} {\rm \cite[Satz]{Wie2}}
Suppose that $\pi$  is a union of  disjoint subsets $\sigma$ and $\tau$.
Assume that a  group $G$ possesses a $\pi$-Hall subgroup $H=H_\sigma\times H_\tau$, where  $H_\sigma$ is a nilpotent
$\sigma$-subgroup and  $H_\tau$ is a $\tau$-subgroup of $H$, and let $G$ satisfy $\D_\tau$. Then $G$ satisfies~$\D_\pi$.
\end{Theo}

This result also is included in the textbooks \cite{Suz,ShemBook,Robinson} and  in the talk  \cite{WieHuppert} and has many generalizations
and analogs (see \cite{Mikaeljan,Rusakov1,Shem5,Shem3,Shem6}).

In the same paper  \cite{Wie2}, Wielandt asked whether, instead of the nilpotency of $H_\sigma$, it is sufficient to assume that $G$ satisfies $\D_\sigma$?
Thus, Wielandt formulated the following conjecture.

\bigskip
\noindent {\bf Conjecture.} (Wielandt,  \cite{Wie2})
{\rm Suppose that a set of primes  $\pi$  is a union of  disjoint subsets $\sigma$ and $\tau$, and a finite group $G$ possesses a $\pi$-Hall subgroup $H=H_\sigma\times H_\tau$, where $H_\sigma$ and $H_\tau$ are
$\sigma$- and   $\tau$-subgroups of $H$, respectively. If $G$ satisfies both $\D_\sigma$ and $\D_\tau$, then $G$ satisfies~$\D_\pi$.}
\bigskip

This conjecture was mentioned in \cite{Hart2,Shem5,Shem3,Shem6} and  was investigated in \cite{Shem5,Shem3,Shem6}. 
The main goal of the present paper is to prove the following theorem which, in particular,  completely confirms Wielandt's conjecture.

\begin{Theo}\label{main1} {\em (Main Theorem)}
Let a set $\pi$ of primes be a union of disjoint subsets $\sigma$ and~$\tau$. Assume that a finite group  $G$ possesses a  $\pi$-Hall subgroup $H=H_\sigma\times H_\tau$, where $H_\sigma$ and $H_\tau$ are
$\sigma$- and   $\tau$-subgroups, respectively. Then  $G$ satisfies~$\D_\pi$ if and only if $G$ satisfies both $\D_\sigma$ and~$\D_\tau$.
\end{Theo}

 In view of Theorem~\ref{W1954}, one can regard Theorem~\ref{main1} as generalization of Theorem  \ref{W1959}.
Moreover, by induction, it is easy to show that  Theorem \ref{main1} is equivalent to the following statement.

\begin{Cor}\label{Cor1}
{ Suppose a set $\pi$ of primes is a union of pairwise disjoint subsets $\pi_i$, $i=1,\dots,n$.
Assume a finite group $G$ possesses a  $\pi$-Hall subgroup 
\begin{linenomath}
$$H=H_{1}\times\dots\times H_{n}$$
\end{linenomath}
where $H_{i}$ is a
$\pi_i$-subgroup  for every  $i=1,\ldots,n$.  Then $G$ satisfies~$\D_\pi$ if and only if  $G$ satisfies  $\D_{\pi_i}$  for every  $i=1,\dots,n$.}
\end{Cor}

In view of Sylow's Theorem, Corollary~\ref{Cor1} is a generalization of  Theorem \ref{W1954} since the statement of the theorem can be obtained if we assume that $\pi_i=\{p_i\}$ for $i=1,\ldots,n$, where 
\begin{linenomath}
$$\{p_1,\ldots,p_n\}=\pi\cap \pi(G).$$
\end{linenomath}

\section{Notation and preliminary results}

Our notation is standard. Given a group $G$ we denote the set of $\pi$-Hall subgroups of $G$ by $\Hall_\pi(G)$. According to \cite{Hall}, the class of  groups $G$ with $\Hall_\pi(G)\ne\varnothing$ is denoted by~$\E_\pi$. We denote by $\Sym_n$ and $\Alt_n$ the symmetric and the alternating groups of degree $n$, respectively. For sporadic groups we use the notation from \cite{Atlas}. The notation for groups of Lie type agrees with that in \cite{Car}. A finite field of cardinality $q$ is denoted by~$\F_q$.

Let $r$ be an odd prime and $q$ be an integer coprime to~$r$. We denote by $e(q,r)$ the minimal natural $e$ such that 
\begin{linenomath}
$${q^e\equiv 1\pmod r},$$
\end{linenomath}
that is, $e(q,r)$ is the multiplicative order of $q$ modulo~$r$.

\begin{Lemma}\label{HallSubgroup} {\em \cite[Lemma~1]{Hall}}
Let $A$ be a normal subgroup and let $H$ be a $\pi$-Hall subgroup of a group $G$. Then
\begin{linenomath}
$${H\cap A\in \Hall_\pi(A)}\text{ and } HA/A\in \Hall_\pi(G/A).$$
\end{linenomath}
\end{Lemma}

\begin{Lemma}\label{DpiExt} {\em \cite[Theorem~7.7]{VR4}} Let $A$ be a normal subgroup of  $G$. Then  
\begin{linenomath}$$G\in\D_\pi\text{ if and only if } A\in \D_\pi \text{ and } G/A\in \D_\pi.$$
\end{linenomath}
\end{Lemma}

We need the criterion obtained in \cite{R4}  for a simple  group  satisfying~$\D_\pi$. In order to formulate it, we need, according with \cite{R4}, to define Conditions I--VII for a pair $(G,\pi)$, where $G$ is a finite simple group and $\pi$ is a set of primes.

\medskip

\noindent{\bf Condition I.}  We say that $(G,\pi)$ satisfies Condition I
if
\begin{linenomath}
$$\text{ either } \pi(G)\subseteq\pi\text{ or }|\pi\cap\pi(G)|\leq 1.$$
\end{linenomath}

\medskip

\noindent{\bf Condition II.} We say that $(G,\pi)$  satisfies
Condition  II if one of the following cases holds.
\begin{itemize}
\item[$(1)$]  $G\simeq M_{11}$ and $\pi\cap\pi(G)=\{5,11\}$;
\item[$(2)$]  $G\simeq M_{12}$  and $\pi\cap\pi(G)=\{5,11\}$;
\item[$(3)$]  $G\simeq M_{22}$ and $\pi\cap\pi(G)=\{5,11\}$;
\item[$(4)$]  $G\simeq M_{23}$  and $\pi\cap\pi(G)$
coincide with one of the following sets $\{5,11\}$ and $\{11,23\}$;
\item[$(5)$]  $G\simeq M_{24}$  and $\pi\cap\pi(G)$
coincide with one of the following sets $\{5,11\}$ and $\{11,23\}$;
\item[$(6)$]  $G\simeq J_1$ and $\pi\cap\pi(G)$ coincide with one of the following sets 
\begin{linenomath}
$$\{3,5\},\,\,
\{3,7\},\,\,  \{3,19\},\,\, \text{ and }\,\, \{5,11\};$$
\end{linenomath}
\item[$(7)$]  $G\simeq J_4$ and $\pi\cap\pi(G)$ coincide with one of the following sets  
\begin{linenomath}
$$\{5,7\},\,\,
\{5,11\},\,\, \{5,31\},\,\, \{7,29\},\,\,\text{ and }\,\, \{7,43\};$$
\end{linenomath}
\item[$(8)$]  $G\simeq O'N$ and $\pi\cap\pi(G)$
coincide with one of the following sets $\{5,11\}$ and $\{5,31\}$;
\item[$(9)$]  $G\simeq Ly$ and $\pi\cap\pi(G)=\{11,67\}$;
\item[$(10)$]  $G\simeq Ru$ and $\pi\cap\pi(G)=\{7,29\}$;
\item[$(11)$]  $G\simeq Co_1$ and $\pi\cap\pi(G)=\{11,23\}$;
\item[$(12)$]  $G\simeq Co_2$  and $\pi\cap\pi(G)=\{11,23\}$;
\item[$(13)$]  $G\simeq Co_3$ and $\pi\cap\pi(G)=\{11,23\}$;
\item[$(14)$]  $G\simeq M(23)$  and $\pi\cap\pi(G)=\{11,23\}$;
\item[$(15)$]  $G\simeq M(24)'$  and $\pi\cap\pi(G)=\{11,23\}$;
\item[$(16)$]  $G\simeq B$  and $\pi\cap\pi(G)$ coincide with one of the following sets $\{11,23\}$ and  $\{23,47\}$;
\item[$(17)$]  $G\simeq M$ and $\pi\cap\pi(G)$ coincide with one of the following sets $\{23,47\}$ and  $\{29,59\}$.
\end{itemize}

\medskip

\noindent{\bf Condition III.} Let $G$ be isomorphic to a group of Lie type over the field $\F_q$ of characteristic $p\in\pi$ and let
\begin{linenomath}
$$\tau=(\pi\cap\pi(S))\setminus\{p\}.$$
\end{linenomath}
We say that $(G,\pi)$
satisfies Condition~III if $\tau\subseteq\pi(q-1)$ and every
prime in $\pi$ does not divide the order of the Weyl group of~$G$.

\medskip
\noindent{\bf Condition IV.} Let $G$ be  isomorphic to   a group of Lie type
with the base field $\F_q$ of  characteristic $p$ but not isomorphic to  ${^2B}_2(q)$, ${^2F}_4(q)$ and ${^2G}_2(q)$. Let
$2,p\not\in\pi$. Denote by $r$ the minimum in $\pi\cap\pi(G)$ and let 
\begin{linenomath}
$$\tau=(\pi\cap\pi(G))\setminus\{r\}\,\, \text{ and  }\,\, a=e(q,r).$$
\end{linenomath}
We say that
$(G,\pi)$  satisfies Condition IV if there exists $t\in\tau$
with $b=e(q,t)\ne a$ and one of the following statements holds.

\begin{itemize}

\item[$(1)$]\noindent  $G\simeq A_{n-1}(q)$, $a=r-1$, $b=r$, $(q^{r-1}-1)_r=r$,
$\left[\displaystyle\frac{n}{r-1}\right]=\left[\displaystyle\frac{n}{r}\right]$,
and $e(q,s)=b$ for every $s\in\tau$;

\item[$(2)$]  $G\simeq A_{n-1}(q)$, $a=r-1$, $b=r$, $(q^{r-1}-1)_r=r$,
$\left[\displaystyle\frac{n}{r-1}\right]=\left[\displaystyle\frac{n}{r}\right]+1$,
$n\equiv -1 \pmod r$, and $e(q,s)=b$ for every $s\in\tau$;

\item[$(3)$]  $G\simeq {^2A}_{n-1}(q)$, $r\equiv 1 \pmod 4$, $a=r-1$, $b=2r$, $(q^{r-1}-1)_r=r$,
$\left[\displaystyle\frac{n}{r-1}\right]=\left[\displaystyle\frac{n}{r}\right]$
and $e(q,s)=b$ for every $s\in\tau$;

\item[$(4)$]  $G\simeq {^2A}_{n-1}(q)$, $r\equiv 3 \pmod 4$, $a=\displaystyle\frac{r-1}{2}$,
$b=2r$, $(q^{r-1}-1)_r=r$,
$\left[\displaystyle\frac{n}{r-1}\right]=\left[\displaystyle\frac{n}{r}\right]$
and  $e(q,s)=b$ for every $s\in\tau$;

\item[$(5)$]  $G\simeq {^2A}_{n-1}(q)$, $r\equiv 1 \pmod 4$, $a=r-1$, $b=2r$, $(q^{r-1}-1)_r=r$,
$\left[\displaystyle\frac{n}{r-1}\right]=\left[\displaystyle\frac{n}{r}\right]+1$,
$n\equiv -1 \pmod r$ and $e(q,s)=b$ for every $s\in\tau$;

\item[$(6)$]  $G\simeq {^2A}_{n-1}(q)$, $r\equiv 3 \pmod 4$, $a=\displaystyle\frac{r-1}{2}$,
$b=2r$, $(q^{r-1}-1)_r=r$,
$\left[\displaystyle\frac{n}{r-1}\right]=\left[\displaystyle\frac{n}{r}\right]+1$,
$n\equiv -1 \pmod r$ and
$e(q,s)=b$ for every  $s\in\tau$;

\item[$(7)$]  $G\simeq {^2D}_n(q)$, $a\equiv 1 \pmod 2$, $n=b=2a$ and for every $s\in\tau$ either $e(q,s)=a$ or
$e(q,s)=b$;

\item[$(8)$]  $G\simeq {^2D}_n(q)$, $b\equiv 1 \pmod 2$, $n=a=2b$ and for every $s\in\tau$ either $e(q,s)=a$ or
$e(q,s)=b$.
\end{itemize}

\medskip

\noindent{\bf Condition V.} Let $G$ be  isomorphic to  a group of Lie type
with the base field $\F_q$ of characteristic $p$, and not isomorphic to ${^2B}_2(q)$, ${^2F}_4(q)$ and ${^2G}_2(q)$. Suppose,
$2,p\not\in\pi$.  Let $r$ be the minimum in $\pi\cap\pi(G)$, let 
\begin{linenomath}
$$\tau=(\pi\cap\pi(G))\setminus\{r\}\,\,\text{ and }\,\, c=e(q,r).$$
\end{linenomath} 
We say that
$(G,\pi)$ satisfies Condition V if $e(q,t)=c$ for every
$t\in\tau$ and one of the following statements holds.
\begin{itemize}
\item[$(1)$]  $G\simeq A_{n-1}(q)$ and $n<cs$ for every $s\in\tau$;

\item[$(2)$]  $G\simeq{^2A}_{n-1}(q)$, $c\equiv 0 \pmod 4$ and $n<cs$ for every $s\in\tau$;

\item[$(3)$]  $G\simeq{^2A}_{n-1}(q)$, $c\equiv 2 \pmod 4$ and ${2n<cs}$ for every $s\in\tau$;

\item[$(4)$]  $G\simeq{^2A}_{n-1}(q)$, $c\equiv 1 \pmod 2$ and ${n<2cs}$ for every $s\in\tau$;

\item[$(5)$]  $G$ is isomorphic to one of the groups $B_n(q)$, $C_n(q)$, or ${^2D}_n(q)$,  $c$
is odd and $2n<cs$ for every $s\in\tau$;

\item[$(6)$]   $G$ is isomorphic to one of the groups $B_n(q)$, $C_n(q)$, or $D_n(q)$, $c$
is even and $n<cs$ for every $s\in\tau$;

\item[$(7)$]  $G\simeq D_n(q)$,  $c$ is even and  $2n\leq cs$ for every $s\in\tau$;

\item[$(8)$]  $G\simeq {^2D}_n(q)$,  $c$  is odd and $n\leq cs$ for every $s\in\tau$;

\item[$(9)$]  $G\simeq {^3D}_4(q)$;

\item[$(10)$]  $G\simeq E_6(q)$,  and if $r=3$ and $c=1$ then $5,13\not\in\tau$;

\item[$(11)$]  $G\simeq {^2E}_6(q)$, and if $r=3$ and  $c=2$ then $5,13\not\in\tau$;

\item[$(12)$]  $G\simeq E_7(q)$, if $r=3$ and $c\in\{1,2\}$ then $5,7,13\not\in\tau$,
and if $r=5$ and $c\in\{1,2\}$ then
$7\not\in\tau$;

\item[$(13)$]  $G\simeq E_8(q)$, if
 $r=3$ and $c\in\{1,2\}$ then $5,7,13\not\in\tau$, and if $r=5$ and
$c\in\{1,2\}$ then $7,31\not\in\tau$;

\item[$(14)$]  $G\simeq G_2(q)$;

\item[$(15)$]  $G\simeq F_4(q)$, and if $r=3$ and $c=1$ then $13\not\in\tau$.
\end{itemize}

\medskip

\noindent{\bf Condition VI.} We say that $(G,\pi)$ satisfies Condition VI
if one of the following statements holds.

\begin{itemize}
\item[$(1)$]  $G$ is isomorphic to ${^2B}_2(2^{2m+1})$ and  $\pi\cap\pi(G)$ is contained in one of the sets
\begin{linenomath}
$$\pi(2^{2m+1}-1), \,\,\, \pi(2^{2m+1}\pm 2^{m+1}+1);$$
\end{linenomath}

\item[$(2)$]  $G$ is isomorphic to ${^2G}_2(3^{2m+1})$ and  $\pi\cap\pi(G)$ is contained in one of the sets
\begin{linenomath}
$$\pi(3^{2m+1}-1)\setminus\{2\}, \,\,\, \pi(3^{2m+1}\pm 3^{m+1}+1)\setminus\{2\};$$
\end{linenomath}

\item[$(3)$]  $G$ is isomorphic to ${^2F}_4(2^{2m+1})$ and $\pi\cap\pi(S)$ is contained in one of the sets
\begin{linenomath}
$$\pi(2^{2(2m+1)}\pm 1), \,\,\, \pi(2^{2m+1}\pm 2^{m+1}+1),$$
\end{linenomath}
\begin{linenomath}
$$\pi(2^{2(2m+1)}\pm 2^{3m+2}\mp2^{m+1}-1),\,\,\, \pi(2^{2(2m+1)}\pm
2^{3m+2}+2^{2m+1}\pm2^{m+1}-1).$$
\end{linenomath}
\end{itemize}

\medskip

\noindent{\bf Condition VII.} Let $G$ be  isomorphic to  a group of Lie type with the base field
$\F_q$ of characteristic $p$. Suppose that  $2\in \pi$ and  $3, p\not\in\pi$, and let
\begin{linenomath} $$\tau=(\pi\cap\pi(G))\setminus\{
2\}\,\,\text{ and }\,\, \varphi=\{t\in\tau\mid t \text{ is a Fermat number}\}.$$
\end{linenomath}
We say
that $(G,\pi)$  satisfies Condition VII if
\begin{linenomath}
$$\tau\subseteq\pi(q-\varepsilon),$$
\end{linenomath}
 where the number $\varepsilon=\pm 1$ is  such that $4$ divides $q-\varepsilon$, and one of the following statements
holds.
\begin{itemize}
\item[$(1)$] $G$ is isomorphic to either ${A}_{n-1}(q)$ or  ${}^2{A}_{n-1}(q)$, $s>n$
 for every $s\in\tau$, and $t>n+1$ for every $t\in\varphi$;
\item[$(2)$] $G\simeq {B}_n(q)$, and $s>2n+1$
for every $s\in\tau$;
\item[$(3)$] $G\simeq {C}_n(q)$, $s>n$ for every
  $s\in\tau$, and $t>2n+1$ for every  $t\in\varphi$;
\item[$(4)$] $G$ is isomorphic to either ${D}_n(q)$ or  ${}^2{D}_n(q)$,  and $s>2n$
for every $s\in\tau$;
\item[$(5)$] $G$ is isomorphic to either ${G}_2(q)$ or ${^2G}_2(q)$, and $7\not\in\tau$;
\item[$(6)$] $G\simeq {F}_4(q)$ and $5,7\not\in\tau$;
\item[$(7)$] $G$ is isomorphic to either $E_6(q)$ or ${}^2E_6(q)$, and $5,7\not\in\tau$;
\item[$(8)$] $G\simeq {E}_7(q)$ and $5,7,11\not\in\tau$;
\item[$(9)$] $G\simeq {E}_8(q)$ and $5,7,11,13\not\in\tau$;
\item[$(10)$] $G\simeq {^3D}_4(q)$ and $7\not\in\tau$.
\end{itemize}

\begin{Lemma}\label{Dpi_ariLemma} {\em \cite[Theorem 3]{R4}} Let $\pi$ be a set of primes and $G$ be a simple group. Then $G\in \D_\pi$ if and only if $(G,\pi)$ satisfies one of Conditions~{\rm I--VII}.
\end{Lemma}

\begin{Lemma}\label{Dpi23}
Let  $G$ be a simple group and  $\pi$ be  such that $2,3\in\pi\cap \pi(G).$ Then 
\begin{linenomath}
$$G\in \D_\pi \text{ if and only if }\pi(G)\subseteq\pi.$$
\end{linenomath}
\end{Lemma}

\begin{proof} If  $\pi(G)\subseteq \pi$, then evidently $G\in \D_\pi$. Conversely, suppose  $G\in \D_\pi$. Then Lemma~\ref{Dpi_ariLemma} implies that  $(G,\pi)$ satisfies one of Conditions I--VII. Without loss of generality, we may assume that  $\pi\cap\pi(G)=\pi$.
\smallskip

If Condition I holds, then 
\begin{linenomath}$$  \text{ either } |\pi|\leq 1 \text{ or } \pi=\pi(G).$$
\end{linenomath}
 Since  $2,3\in\pi$, we have that $|\pi|\geq 2$ and thus $\pi=\pi(G)$.\smallskip

Clearly, $(G,\pi)$ cannot satisfy one of Conditions II and IV--VII since otherwise either $2\notin\pi$ (Conditions II and IV--VI) or $3\notin\pi$ (Condition~VII).\smallskip

Finally, $(G,\pi)$ also cannot satisfy Condition III. Indeed, if Condition III holds, then $G$ is a group of Lie type over a field of characteristic $p\in\pi$ and every prime in $\pi$ (in particular $2$) does not divide the order of the Weyl group of $G$. But this is impossible since the Weyl group is nontrivial and is generated by involutions (see~\cite[page~13 and Proposition~13.1.2]{Car}).
\end{proof}

\begin{Lemma}\label{HallSymmetric}
Suppose that $n\geq 5$ and $\pi$ is a set of primes with 
\begin{linenomath}
$$|\pi\cap \pi(n!)|>1 \text{ and } \pi(n!)\nsubseteq\pi.$$
\end{linenomath}
Then
\begin{itemize}
 \item[{\em (1)}] The complete list of possibilities for $\Sym_n$  containing a
$\pi$-Hall subgroup $H$ is given in Table~{\em\ref{Symmetric}}.
 \item[{\em (2)}] $K\in\Hall_\pi(\Alt_n)$ if and only if  ${K=H\cap
\Alt_n}$ for some ${H\in\Hall_\pi(\Sym_n)}$.
\end{itemize}
In particular, if either $2\notin \pi$ or $3\notin\pi$, then $\Alt_n\notin\E_\pi$.
\end{Lemma}

\begin{table}\caption{\label{Symmetric}}
\begin{center}
\begin{tabular}{|c|c|c|}
\hline
$n$&$\pi$&$H\in\Hall_\pi(\Sym_n)$\\ \hline\hline
prime&$\pi((n-1)!)$&$\Sym_{n-1}$\\
$7$&$\{2,3\}$&$\Sym_3\times\Sym_4$\\
$8$&$\{2,3\}$&$\Sym_4\wr\Sym_2$\\ \hline
\end{tabular}
\end{center}
\end{table}
\begin{proof} See  \cite[Theorem~A4 and the notices after it]{Hall}, \cite[Main result]{Thom}, and \cite[Theorem~4.3 and Corollary~4.4]{VR4}.\end{proof}

\begin{Lemma}
\label{HallSporadic} {\em \cite[Corollary 6.13]{Gross1}}
Let $G$~be either one of~$26$ sporadic groups or the Tits group. Assume that $2\notin\pi$. Then $G\in \E_\pi$ if and only if either $|\pi\cap \pi(G)|\leq 1$ or $G$ and $\pi\cap\pi(G)$ are given in Table~{\em
\ref{Sporadicodd}}. In particular, $|\pi\cap \pi(G)|\leq 2$.
\end{Lemma}

\begin{table}\caption{}\label{Sporadicodd}
\begin{center}
\begin{tabular}{lcr}%
\begin{tabular}{|c|c|}
\hline
 $G$& $\pi\cap\pi(G)$ \\ \hline \hline
 $M_{11}$& $\{5,11\}$  \\ \hline
 $M_{23}$&$\{5,11\}$ \\
 &$\{11,23\}$\\ \hline
$Ru$&$\{7,29\}$\\
\hline
$Fi_{23}$& $\{11,23\}$\\ \hline
$J_1$&$\{3,5\}$\\
&$\{3,7\}$\\
& $\{3,19\}$\\
&$\{5,11\}$\\
\hline
$Co_1$&$\{11,23\}$\\ \hline

\end{tabular}\,\,\,\,\,
\begin{tabular}{|c|c|}
\hline
$G$ &$\pi\cap\pi(G)$ \\
\hline\hline
$M_{12}$&$\{5,11\}$\\ \hline
 $M_{24}$&$\{5,11\}$\\
&$\{11,23\}$\\ \hline
 $Fi_{24}'$&$\{11,23\}$  \\ \hline
$J_4$&$\{5,7\}$ \\
&$\{5,11\}$\\
&$\{5,31\}$\\
&$\{7,29\}$\\
&$\{7,43\}$ \\\hline
$Co_2$&$\{11,23\}$\\ \hline

\end{tabular}\,\,\,\,\,
\begin{tabular}{|c|c|}
\hline
 $G$ & $\pi\cap\pi(G)$\\ \hline\hline
 $M_{22}$& $\{5,11\}$\\ \hline
  $ Ly$& $\{11,67\}$\\ \hline
 $O'N$&$\{3,5\}$\\
 &$\{5,11\}$ \\
&$\{5,31\}$ \\
\hline

$B$&$\{11,23\}$\\
&$\{23,47\}$\\ \hline
$M$&$\{23,47\}$\\
&$\{29,59\}$ \\ \hline
$Co_3$&$\{11,23\}$\\
\hline
\end{tabular}
\end{tabular}
\end{center}
\end{table}

\begin{Lemma}\label{HallSpor2} {\em \cite[Theorem 4.1]{R}}
Let $G$ be either a simple sporadic group or the Tits group and   $\pi$ be  such that 
\begin{linenomath}
$$2\in\pi,\,\, \pi(G)\nsubseteq\pi, \text{ and } |\pi\cap\pi(G)|>1.$$ 
\end{linenomath}
Then
the complete list for $G$  containing a $\pi$-Hall subgroup $H$ is given in Table~{\em\ref{tb0}}.
In particular, if   $3\notin\pi$, then  
\begin{linenomath}
$$G= J_1 \text{ and } \pi\cap\pi(G)=\{2,7\}.$$
\end{linenomath}
\end{Lemma}

\begin{longtable}{|l|l|r|}\caption{}\label{tb0}\\  \hline
$G$ & $\pi\cap\pi(G)$ &
Structure of $H$ \\ \hline\hline $M_{11}$  & $\{2,3\}$ & $3^2\splitext Q_8\arbitraryext 2$ \\ &
$\{2,3,5\}$       &  $\Alt_6\arbitraryext 2$\\ \hline $M_{22}$  & $\{2,3,5\}$ &
$2^4\splitext \Alt_6$\\ \hline $M_{23}$  & $\{2,3\}$         & $2^4\splitext(3\times
A_4)\splitext2$\\ & $\{2,3,5\}$ &  $2^4\splitext \Alt_6$\\ & $\{2,3,5\}$       &
$2^4\splitext (3\times \Alt_5)\splitext 2$\\ & $\{2,3,5,7\}$     &  ${\rm
          L}_3(4)\splitext 2_2$\\ & $\{2,3,5,7\}$     &  $2^4\splitext \Alt_7$\\ &
          $\{2,3,5,7,11\}$  & $M_{22}$\\ \hline $M_{24}$  & $\{2,3,5\}$
          & $2^6 \splitext 3\nonsplitext \Sym_6$\\ \hline $J_1$     & $\{2,3\}$         &
$2\times \Alt_4$\\ & $\{2,7\}$         &  $2^3\splitext 7$\\ & $\{2,3,5\}$ &
          $2\times A_5$\\ & $\{2,3,7\}$       &  $2^3\splitext 7\splitext 3$\\ \hline
$J_4$     & $\{2,3,5\}$       &  $2^{11}\splitext (2^6\splitext 3\nonsplitext \Sym_6)$\\
\hline
\end{longtable}

\begin{Lemma}
\label{p_ip_pi}
Let  $G$ be a group of Lie type with base field $\F_q$ of characteristic $p$. Assume that  $\pi$  is such that  $p\in\pi$, and either $2\notin\pi$ or $3\notin\pi$. Suppose  $G\in \E_\pi$ and $H\in\Hall_\pi(G)$. Then one of the following statements holds.
\begin{itemize}
\item[$(1)$] $\pi\cap \pi(G)\subseteq\pi(q-1)\cup\{p\}$, a Sylow $p$-subgroup $P$ of $H$ is normal in  $H$ and  $H/P$ is~abelian.
\item[$(2)$] $p=2$,  $G\simeq {}^2B_2(2^{2n+1})$ and $\pi(G)\subseteq\pi$.
\end{itemize}
\end{Lemma}

\begin{proof} See \cite[Theorem 3.2]{Gross1} and \cite[Theorem~3.1]{Gross4}.\end{proof}

\begin{Lemma}
\label{p_notip_pi}
Let  $G$ be a group of Lie type over a field
of characteristic  $p$.  Assume that  $\pi$  is such that $p\notin\pi$, and either  $2\notin\pi$ or $3\notin\pi$. Denote by  $r$ the minimum in $\pi\cap \pi(G)$. Suppose $G\in \E_\pi$ and $H\in\Hall_\pi(G)$. Then $H$ possesses a normal abelian $r'$-Hall subgroup.
\end{Lemma}

\begin{proof} See \cite[Theorems 4.6 and 4.8, and Corollary  4.7]{Gross2}, \cite[Theorem~1]{VR2}, and \cite[Lem\-ma~5.1 and Theorem~5.2]{VR1}.\end{proof}

\begin{Lemma}
\label{2or3notin_pi}
Let  $G$ be a simple nonabelian group. Assume that  $\pi$  is such that $\pi(G)\not\subseteq\pi$ and either $2\notin\pi$ or $3\notin\pi$. Suppose $G\in \E_\pi$ and $H\in\Hall_\pi(G)$. Then  $H$ is solvable and, for any partition 
\begin{linenomath}
$$\pi\cap\pi(G)=\sigma\cup\tau,$$ 
\end{linenomath}
where $\sigma$ and $\tau$ are disjoint nonempty sets, either  $\sigma$-Hall or $\tau$-Hall subgroup of  $H$ is nilpotent.
\end{Lemma}

\begin{proof}
Consider all possibilities, according to the classification of finite simple groups (see \cite[Theorem~0.1.1]{AschLyoSmSol}).\smallskip

C a s e~~1:~~$G\simeq \Alt_n$, $n\geq 5$. By Lemma~\ref{HallSymmetric} it follows that $|\pi\cap \pi(G)|=1$ and a partition  $\pi\cap\pi(G)=\sigma\cup\tau$ onto nontrivial disjoint subsets is impossible.\smallskip

C a s e~~ 2:~~$G$ is either a sporadic group or the Tits group. By Lemmas  \ref{HallSporadic} and \ref{HallSpor2} it follows that either  
\begin{linenomath}
$$|\pi\cap \pi(G)|=1,$$ or $$2\notin\pi \text { and } |\pi\cap \pi(G)|=2,$$ 
\end{linenomath} 
or 
\begin{linenomath}
$$3\notin\pi, \,\, G\simeq J_1 \text{ and } \pi\cap\pi(G)=\{2,7\}.$$ 
\end{linenomath}
If  $|\pi\cap\pi(G)|=1$, then a partition  
\begin{linenomath}
$$\pi\cap \pi(G)=\sigma\cup\tau$$ 
\end{linenomath}
onto nonempty  disjoint subsets is impossible. If $|\pi\cap \pi(G)|=2$, then $H$ is solvable by Burnside's $p^aq^b$-theorem \cite[Ch.~I,~2]{DH}, the orders of its $\sigma$-Hall and  $\tau$-Hall subgroups are  powers of  primes, and thus  $\sigma$-Hall and $\tau$-Hall subgroups of $G$ are nilpotent.\smallskip

C a s e~~ 3:~~$G$ is a group of Lie type over a field of characteristic $p\in\pi$. We may assume, without loss of generality, that $|\pi\cap \pi(G)|>1$ and $p\in\sigma$. By Lemma \ref{p_ip_pi}, $H$ is solvable and its  $\tau$-Hall subgroup $T$ is isomorphic to a subgroup of the abelian group $H/P$, where $P$ is the (normal) Sylow  $p$-subgroup of $H$. Whence   $T$ is abelian and, in particular, is nilpotent.\smallskip

C a s e~~ 4:~~$G$ is a group of Lie type over a field of characteristic $p\notin\pi$. We may assume, without loss of generality, that $|\pi\cap\pi(G)|>1$. Denote by $r$
the minimum in $\pi\cap\pi(G)$, and assume that $r\in\sigma$. By Lemma \ref{p_notip_pi}, it follows that $H$ is solvable and its $\tau$-Hall subgroup $T$ is included in the normal abelian  $r'$-Hall subgroup of $H$. Thus we again obtain that  $T$ is abelian.\smallskip

Thus in the all cases, Lemma~\ref{2or3notin_pi} holds.
\end{proof}

\section{Proof of the main theorem}

Assume that the hypothesis of Theorem \ref{main1} holds, i.e. we have a partition  
\begin{linenomath}
$$\pi=\sigma\cup\tau$$
\end{linenomath} 
of $\pi$ onto disjoint subsets  $\sigma$ and $\tau$, and a group  $G$ satisfying condition:
\begin{itemize}
\item[$(1)$] $G$ possesses a $\pi$-Hall subgroup  $H$ such that  
\begin{linenomath}
$$H=H_\sigma\times H_\tau,$$ 
\end{linenomath}
where  $H_\sigma$ and $H_\tau$ are $\sigma$- and $\tau$-subgroups, respectively.
\end{itemize}
It is easy to see that  $H_\sigma$ and $H_\tau$ are, respectively, $\sigma$-Hall and $\tau$-Hall subgroups of both $H$ and $G$. Moreover, $H_\sigma$ coincides with the set of all $\sigma$-elements of $H$, while $H_\tau$ is the set of all $\tau$-elements of~$H$.

We prove first that $G\in\D_\pi$ implies $G\in \D_\sigma\cap\D_\tau$. We need to prove that a $\sigma$-subgroup $K$ of $G$ is conjugate to a subgroup of $H_\sigma$ in order to prove that $G\in \D_\sigma$. Since $K$ is, in particular, a $\pi$-subgroup and $G\in \D_\pi$, there exists $g\in G$ such that $K^g\leq H$. Hence $K^g\leq H_\sigma$, since $H_\sigma$ is the set of all $\sigma$-elements of $H$. Thus we obtain $G\in \D_\sigma$.  The same arguments show that $G\in D_\tau$.

Now we prove the converse statement: if $G\in\D_\sigma\cap \D_\tau$, then $G\in\D_\pi$. Assume that it fails. Without loss of generality we may assume that $G$ satisfies the following conditions:
\begin{itemize}
\item[$(2)$]  $G\in\D_\sigma\cap\D_\tau$;
\item[$(3)$] $G\notin\D_\pi$;
\item[$(4)$] $G$ has the smallest possible order in the class of groups satisfying conditions $(1)$--$(3)$.
\end{itemize}
Now we show that the assumption of existence of such group leads us to a contradiction.

In view of (3) we have $\pi(G)\nsubseteq\pi$ and $G$ is nonabelian.

Note that $G$ must be simple. Indeed, assume that  $G$ possesses a nontrivial proper normal subgroup  $A$. Then Lemma \ref{HallSubgroup} impies that $A$ and $G/A$ satisfy $(1)$, and Lemma  \ref{DpiExt} implies that they both satisfy $(2)$. In view of $(4)$, neither  $A$ nor $G/A$ satisfies  $(3)$, and thus  $A\in\D_\pi$ and $G/A\in\D_\pi$. Hence by Lemma~\ref{DpiExt}, we obtain $G\in\D_\pi$, which contradicts the assumption~$(3)$.

Assume that either $2$ or $3$ does not lie in $\pi$. Then by Lemma \ref{2or3notin_pi} either  $H_\sigma$ or $H_\tau$ is nilpotent. Hence by Theorem~\ref{W1959} we obtain $G\in\D_\pi$, which contradicts~(3). Hence $2,3\in\pi$.

By Lemma  \ref{Dpi23} and the condition~$(2)$,  the numbers $2$ and $3$ cannot simultaneously lie in the same subset  $\sigma$ or $\tau$. We may, therefore, assume that  \begin{linenomath}
$$2\in\sigma \text{ and } 3\in\tau.$$
\end{linenomath}

Let $S$ be a Sylow  2-subgroup of  $H_\sigma$ (hence of both $H$ and $G$), and $T$ be a Sylow  $3$-subgroup of $H_\tau$ (hence of both $H$ and $G$).  Since   
\begin{linenomath}
$$[S,T]\leq [H_\sigma,H_\tau]=1,$$ 
\end{linenomath}
we see that $G$ possesses a nilpotent $\{2,3\}$-Hall subgroup   
\begin{linenomath}$$\langle S, T\rangle\simeq S\times T.$$ 
\end{linenomath}
It follows from Theorem \ref{W1954} that $G\in\D_{\{2,3\}}$. Now by Lemma~\ref{Dpi23} we have $$\pi(G)=\{2,3\}\subseteq\pi,$$ which implies that $G$ is solvable by Burnside's $p^aq^b$-theorem  \cite[Ch.~I,~2]{DH}. This contradiction completes the proof.\qed

\medskip

Notice that the proof of Theorem~\ref{main1} implies the following statement.

\begin{Cor}
Suppose that a set $\pi$ of primes is a disjoint union of subsets $\sigma$ and $\tau$.
Suppose that a finite simple group $G$ possesses a $\pi$-Hall subgroup $H=H_\sigma\times H_\tau$, where $H_\sigma$ and $H_\tau$ are
$\sigma$- and   $\tau$-subgroups, respectively. If $G$ satisfies both $\D_\sigma$ and $\D_\tau$, then either $H_\sigma$ or $H_\tau$ is nilpotent.
\end{Cor}


\end{document}